\documentclass[12pt]{amsart}
\usepackage{amssymb,latexsym}
\usepackage{a4wide}

\newtheorem{theorem}{Theorem}

\newtheorem{HPV}{Theorem HPV \!\!\!\!}

\newtheorem{lemma}{Lemma}

\newtheorem{question}{Question}

\newtheorem*{conjecture}{Conjecture}

\newcommand{\twolineprod}[2]{\prod_{\substack{{\scriptstyle #1}\\
{\scriptstyle #2}}}}

\def\ca {{\mathcal A}}
\def\cd {{\mathcal D}}
\def\ce {{\mathcal E}}

\def\beq {\begin{equation}}
\def\endq {\end{equation}}
\def\g {{\rm{gcd}}}
\def\N {\mathbb{N}}
\def\R {\mathbb{R}}
\def\Z {\mathbb{Z}}
\def\rar {\rightarrow}
\def\E {\mathcal{E}}

\title[Duffin-Schaeffer with extra divergence II]{The Duffin-Schaeffer Conjecture\\ with extra divergence II}

\author[V.~Beresnevich,~G.~Harman,~A.~Haynes,~and~S.~Velani]{Victor Beresnevich, Glyn Harman,\\ Alan Haynes and Sanju Velani}

\begin{document}

\thanks{{\sf Victor Beresnevich:}~University of York, Heslington, York YO10 5DD, United Kingdom.\\ \indent E-mail: {\tt victor.beresnevich@york.ac.uk} \,\,\,\,\,\,\,Tel. +44 1904 32 3072\smallskip}

\thanks{{\sf Glyn Harman:}~Royal Holloway, University of London, Egham, Surrey TW20 0EX, United Kingdom.\\ \indent E-mail: {\tt G.Harman@rhul.ac.uk}\smallskip}

\thanks{{\sf Alan Haynes:}~Department of Mathematics, University Walk, Bristol, BS8 1TW, United Kingdom.\\ \indent E-mail: {\tt alan.haynes@bris.ac.uk}\smallskip}

\thanks{Sanju Velani:~University of York, Heslington, York YO10 5DD, United Kingdom.\\ \indent E-mail: {\tt sanju.velani@york.ac.uk}}

\begin{abstract}
In \cite{HPV} the authors set out a programme to prove the Duffin-Schaeffer Conjecture for measures arbitrarily close to Lebesgue measure. In this paper we take a new step in this direction. Given a nonnegative function $\psi : \N \to \R $, let  $W(\psi)$
denote the set of real numbers $x$ such that  $|nx -a| < \psi(n) $
for infinitely many reduced rationals $a/n \ (n>0) $. Our main result is that $W(\psi)$ is of full Lebesgue measure if there exists a $c > 0 $ such that
$$
\sum_{n\ge 16}  \, \frac{\varphi(n) \psi(n)}{n \exp(c(\log \log n)(\log \log \log n))}  \, =  \, \infty \, .
$$
\\[5ex]
{\it Keywords: Diophantine approximation, Metric number theory, Duffin-Schaeffer conjecture}\\
{\it Mathematics Subject Classification 2000: 11J83, 11K55, 11K60}

\end{abstract}

\maketitle

\section{Introduction}
We use the following notation: $p$ denotes a prime number, $\varphi (n)$ is the Euler phi function, $\lambda$ denotes Lebesgue measure
on $\R/\Z$, and $f\ll g$ means that the absolute value of $f$ is bounded above by a constant times the absolute value of $g$.

Let $\psi:\N\rar\R$ be a non-negative arithmetical function and for
each positive integer $n$ define $\E_n\subseteq \R/\Z$ by
\begin{equation*}
\E_n \, := \,
\bigcup_{\substack{a=1\\\g(a,n)=1}}^n\left(\frac{a-\psi
(n)}{n},\frac{a+\psi (n)}{n}\right).
\end{equation*}
Denote the collection of points $x\in\R/\Z$ which fall in
infinitely many of the  sets $\E_n$ by $W(\psi)$. In other words,
$$
W(\psi):=
\limsup_{n\to\infty}\E_n:=\bigcap_{m=1}^\infty\ \bigcup_{n\ge
m}\E_n   \ .
$$

The question we address is:
\begin{question}
Let $\psi$ be any non-negative arithmetical function.  Under what circumstances is it true that $\lambda (W(\psi))=1$?
\end{question}

It is very easy to give a \emph{necessary} condition for this to happen,
namely the divergence of the series:
\beq\label{4}
\sum_{n=1}^{\infty} \frac{\varphi(n)}{n} \psi(n).
\endq
This follows from the Borel-Cantelli Lemma, since
\[
\lambda(\E_n) \le 2\frac{\varphi(n)}{n} \psi(n)
\]
and so the convergence of (\ref{4})  implies that $\lambda (W(\psi))=0$.
It is a central open problem in metric number theory to show that the divergence of (\ref{4}) is actually \emph{sufficient} to conclude that $\lambda(W(\psi))=1$.
\begin{conjecture}[Duffin-Schaeffer 1941]
We have that $\lambda(W(\psi))=1$ if and only if (\ref{4}) diverges.
\end{conjecture}
There are several significant partial results towards this conjecture, most notably those due to Khintchine, Duffin \& Schaeffer, Erd\"{o}s, Vaaler, and Pollington \& Vaughan \cite{K,Duff,E,V,PV}. The proofs of these results and others are all given in \cite[Chps 2 \& 3]{Harman}.  Recently Pollington and two of this paper's authors \cite{HPV} have considered the effect on the problem of assuming ``extra divergence''. They have posed the following question.

\begin{question}
For what functions $f$ does the divergence of
\beq\label{extra}
\sum_{n=1}^{\infty} f\left(\frac{\psi(n)}{n}\right) \varphi(n)
\endq
guarantee that $\lambda (W(\psi))=1$?
\end{question}

In view of the Mass Transference Principle \cite{MTP} this question is equivalent to investigating the Duffin-Schaeffer Conjecture for (Hausdorff) measures ``arbitrarily'' close to
Lebesgue measure  -- see \cite[\S 5]{HPV} for the details.
Regarding Question 2 itself,  the following result is established in \cite{HPV}.

\begin{HPV}
Let $\psi$ be any non-negative arithmetical function and define the function  $f$ by
\[
f(x) = \begin{cases} 0 &\text{if \ $x = 0$},\\
x \exp\left(\frac{\log x}{\log(-\log x)}  \right) &\text{if \ $0 < x < 1$},\\
1 &\text{if \ $x \ge 1$}.
\end{cases}
\]
Then the divergence of \eqref{extra} is sufficient to conclude that $\lambda(W(\psi))=1$.
\end{HPV}

The authors of \cite{HPV} set as an explicit unsolved problem the task of replacing $f(x)$ above by $f(x) = x(-\log x)^{-1}$.  The hope is that
as one approaches the situation of the Duffin-Schaeffer Conjecture (that is, $f(x) = x$) one can see more clearly
the outstanding problems. In addition any subsequent attack on the conjecture can
assume that the series \eqref{4} is diverging ``slowly" in certain senses.  In this paper we make progress
towards this goal by establishing the following result.

\begin{theorem}
Let $\psi$ be any non-negative arithmetical function and for  any $c > 0$, define the function $f_c$ by
\[
f_c(x) =
\begin{cases} 0 &\text{ \ if \ \  $x = 0$},\\
 x \exp\left(-c (\log (- \log x))(\log \log (-\log x)) \right) &\text{  \ if \ \  $0 < x < 1$},\\
1 &\text{ \ if \ \ $x \ge 1$}.
\end{cases}
\]
Then the divergence of \eqref{extra} with $f = f_c$ is sufficient to conclude that $\lambda(W(\psi))=1$.
\end{theorem}

It is worth pointing out that this strengthening of Theorem HPV is not a consequence of simply  tweaking the approach taken in \cite{HPV} -- see also the remark at the end of \S3 in \cite{HPV}.
%
%
%
By appealing to the
Erd\"{o}s-Vaaler Theorem \cite{V} and to \cite[Theorem
2]{PV} we can assume without loss of generality
throughout the proof that $1/n\le \psi (n)\le 1/2$  whenever $\psi(n) \neq 0$.  In view of this it suffices to prove the following theorem.

\begin{theorem} \label{prove}
The Duffin-Schaeffer Conjecture is true for any non-negative arithmetical function $\psi$ such that the series
\beq\label{nextra}
\sum_{n=16}^{\infty} \frac{\varphi(n) \psi(n)}{n \exp(c(\log \log n)(\log \log \log n))}
\endq
diverges for some $c > 0$.
\end{theorem}

\section{The basic framework}
Gallagher \cite{G} (see also \cite[\S2.2]{Harman}) proved that there is a ``zero-one'' law for Question 1. That is,  for any given
function $\psi$ we have $\lambda(W(\psi))=0$ or $1$.  We therefore only need to prove that under our extra divergence hypothesis $\lambda (W(\psi))>0$. To do this we need the following result \cite[Lemma 2.3]{Harman} whose proof involves little
more than the correct application of the Cauchy-Schwartz inequality.

\begin{lemma}\label{Lemma-quasiind}
Let $\ca_n$ be a sequence of Lebesgue measurable subsets of $\R/\Z$.  Let $\ca$ be the set of $\alpha$ belonging to
infinitely many $\ca_n$.  Then
\beq\label{6}
\lambda(\ca) \ge \limsup_{N \rightarrow \infty} \left(\sum_{n=1}^N \lambda(\ca_n) \right)^2 \left(
\sum_{m,n=1}^N \lambda(\ca_m \cap \ca_m) \right)^{-1}.
\endq
\end{lemma}

The well known Duffin-Schaeffer result \cite[Theorem 2.5]{Harman} toward the Duffin-Schaeffer Conjecture follows from this lemma together with the elementary bound
\beq\label{7}
\lambda(\E_m \cap \E_n) \le 8 \psi(n) \psi(m) ~\text{ for }~ m\not= n.
\endq
However if the sets in the collection $\{\E_n\}$ were quasi-pairwise independent, i.e. if
\beq\label{8}
\lambda(\E_m \cap \E_n) \ll \psi(n) \psi(m) \frac{\varphi(m)}{m} \frac{\varphi(n)}{n} ~\text{ for }~ m\not= n,
\endq
then the Duffin-Schaeffer Conjecture would follow at once from (\ref{6}) together with Gallagher's zero-one law. Unfortunately (\ref{8}) does not hold uniformly for all $m\not= n$. The best known upper bound for $\lambda(\E_n \cap \E_m)$ is the following result.

\begin{lemma}\label{lem2}
For $m \ne n$ we have
\beq\label{9}
\lambda(\E_m \cap \E_n) \ll \lambda(\E_m) \lambda(\E_n)P(m,n),
\endq
where
\beq\label{10}
P(m,n) = \twolineprod{p|mn/\g(m,n)^2}{p >D(m,n)}\left(1 - \frac{1}{p}\right)^{-1},
\endq
with
\beq\label{11}
D(m,n) = \frac{\max(n \psi(m), m \psi(n))}{\g(m,n)}.
\endq
\end{lemma}

This was first stated by Strauch \cite{Strauch}, but was also given independently by Pollington and Vaughan \cite{PV}.
The proof is still essentially elementary, but fairly complicated, and needing a simple sieve upper bound.
Effectively the same result was given earlier by Erd\"os \cite{E}.
Clearly what needs to be done in applying Lemma \ref{Lemma-quasiind} is to show that the factor $P(m,n)$ is bounded
on average.

It is worth pausing here to see what the real difficulties are in estimating $\lambda(\E_m \cap \E_n)$.
Two intervals from $\E_m$ and $\E_n$ overlap if
\[
\left|\frac{a}{m} - \frac{b}{n}  \right|< \frac{\psi(m)}{m} + \frac{\psi(n)}{n}.
\]
We lose nothing in terms of the order of magnitude of the bound in replacing this with
\beq\label{12}
|an-bm| < A(m,n) := 2 \max(m \psi(n), n \psi(m)).
\endq
The length of the intersection is no more than the smallest length
of the two intervals (again nothing is lost in order of magnitude in making this assumption).
We thus have
\[
\lambda(\E_m \cap \E_n)
\le
\ 2 \min\left(\frac{\psi(m)}{m}, \frac{\psi(n)}{n} \right) \Sigma(m,n),
\]
where $\Sigma(m,n)$ denotes the number of solutions to \eqref{12} with
\[
1 \le a < m, \ \ 1 \le b < n, \ \ \g(a,m) = 1, \ \ \g(b,n) = 1.
\]
We thus need to show that, at least on average over $m,n$, we have
\[
\Sigma(m,n) \ll  A(m,n) \frac{\varphi(m)}{m} \frac{\varphi(n)}{n}.
\]
Now if $\g(m,n)=1$ and $A(m,n)$ is not too small there is no problem with this.  The trouble essentially comes when
\[
1 \le \frac{A(m,n)}{\g(m,n)} < \log(mn).
\]
In that case we are not averaging over enough values of $h$ in the equation $an-bm=h$ to
get the required bound.  This is a \emph{real} problem, not just a deficiency in our knowledge.
Our hope would be that
the values of $m$ and $n$ concerned do not make the major contribution to
\[
\sum_{1\le m,n \le N} \lambda(\E_m \cap \E_n).
\]

\section{Proof of Theorem \ref{prove}}
Without loss of generality assume that $\psi (n)\ge n^{-1}$ whenever $\psi (n)\not= 0$  --  see the discussion just before the statement of the theorem.
We divide the integers $n$ into blocks
\[2^{4^h} \le n < 2^{4^{h+1}}.\]
It then follows that
the series \eqref{nextra} diverges with $n$ restricted to blocks with either all the $h$ even, or all the
$h$ odd.  Without loss of generality we suppose the series diverges over blocks with $h$ even, and that $\psi (n)=0$ for all integers $n$ which lie in blocks with $h$ odd.  We then note that
if $m<n$,  if $\psi(m),\psi(n)>0$, and if $m$ and $n$ are in different blocks then
\[
A(m,n) \ge 2n \psi(m) \ge 2n m^{-1} \ge 2n \, \g(m,n)m^{-2} \gg \, \g(m,n) (\log nm).
\]
Hence $P(m,n) \ll 1$ if $m$ and $n$ belong to different blocks.

Now we consider a block $2^{4^h} = X \le m,n < X^4$.  Write $R = \log \log X$ and
\[
\Psi(X) = \sum_{X \le m,n < X^4} \psi(n) \psi(m) \frac{\varphi(m)}{m} \frac{\varphi(n)}{n}.
\]
By one of Mertens' theorems,  we have for $D(m,n) \ge 1$ that
\[
P(m,n) \ll \exp\left(\sum_{D(m,n) < p < \log X} \frac{1}{p}  \right) \ll \frac{R}{1 + \log D(m,n)}.
\]
We let $\cd_j$ be the collection of pairs $(m,n)$ such that $e^j \le D(m,n) < e^{j+1}$. The idea is going to be to divide each $\psi(n)$ by
a suitable factor (say $e^k$) so that the contribution from $R$ consecutive ranges for which we cannot assume $P(m,n)\ll 1$
is not of a larger magnitude than the expected
overall contribution.  We have
\[
\sum_{k \le K} \sum_{k \le j \le k + R} \frac{R}{j+1-k} \sum_{(m,n) \in \cd_j}
\psi(n) \psi(m) \frac{\varphi(m)}{m} \frac{\varphi(n)}{n}
\ll R \Psi(X) \log K,
\]
since each set $\cd_j$ is counted with weight
\[
\le \sum_{k \le \min(j,K)} \frac{R}{j+1-k} \ll R \log K.
\]
We can therefore choose an integer $k \le cR \log R$ such that
\[
\sum_{k \le j \le k +R} \frac{R}{j+1-k} \sum_{(m,n) \in \cd_j}
\psi(n) \psi(m) \frac{\varphi(m)}{m} \frac{\varphi(n)}{n}
\ll \Psi(X).
\]
With this choice of $k$ write
\[
\ce = \bigcup_{k \le j \le k + R} \cd_j.
\]

Now put
\[
\rho(n) = \psi(n)e^{-k} \ \ \text{for} \ X \le n < X^4,
\]
and note that
\[
e^k \le \exp\left(c(\log \log X)(\log \log \log X)\right).
\]

We now assume the above procedure has been carried out on all blocks with $h$ even, and we put $\rho(n)=0$ for
$n$ in a block with odd $h$, and we consider the sets $\E_n(\rho)$. By construction
\[
\sum_{n=1}^{\infty} \frac{\rho(n) \varphi(n)}{n}=\infty,
\]
and
\[
\lambda(\E_m(\rho) \cap \E_n(\rho)) \ll \lambda(\E_m(\rho)) \lambda(\E_n(\rho))
\]
unless $(m,n) \in \ce$. But now we also have that
\[
\sum_{(m,n)\in~\ce} \rho(m)\rho(n) \frac{\varphi(n)}{n} \frac{\varphi(m)}{m} P(m,n) \ll
\sum_{X \le m,n < X^4} \rho(m)\rho(n) \frac{\varphi(n)}{n} \frac{\varphi(m)}{m}
\]
(we note that $P(m,n)$ here does depend on $\rho$) and so
\[
\sum_{m,n=1}^N \lambda(\E_m(\rho) \cap \E_m(\rho)) \ll
\left(\sum_{n=1}^N \lambda(\E_n (\rho)) \right)^2
\]
for $N$ taking the values $2^{4^{h+1}}$.
By Lemma \ref{Lemma-quasiind} and Gallagher's zero-one law we have that $\lambda (W(\rho))=1$, and since $W(\rho)\subseteq W(\psi)$ the proof is completed.

\vspace{6ex}

\noindent{\em Acknowledgements. } 
AH was supported by EPSRC's grant EP/F027028/1. SV was supported by by EPSRC's grants EP/E061613/1 and EP/F027028/1. SV would like to thank Andy Pollington for the many  discussions centered around the Duffin-Schaeffer Conjecture. Also a great thanks to Emma Robertson and Kevin Hall for coaching the girls U9's Wigton Moor football team  to the West Riding girls league and cup double and for putting up with the dynamo duo Ayesha and Iona. I am truly impressed by your dedication, your fairness and above all your friendship and openness with the team and parents.  Thankyou!

\end{document}